\documentclass[12pt]{article}

\usepackage{amssymb, theorem}

\newcommand{\Aalpha}{A_{\alpha}}
\newcommand{\Ae}{A^e}
\newcommand{\Aemod}{\mbox{\bf $\Ae$-mod}}
\newcommand{\Amod}{\mbox{$A$-{\bf mod}}}
\newcommand{\Amodrat}{\mbox{$A$-\bf mod$_\mathbf{rat}$}}
\newcommand{\Balpha}{B^{\alpha}}
\newcommand{\Calpha}{C^{\alpha}}
\newcommand{\Cctn}{C_{\mathrm{ctn}}}
 
\newcommand{\comodDA}{\mbox{{\bf comod}-DA}}
\newcommand{\cotor}{\mathrm{Cotor}}
\newcommand{\enmorph}{\mathrm{End}}
\newcommand{\eg}{{\it e.g.\/}\ }
\newcommand{\Galpha}{G_{\alpha}}
\newcommand{\Hctn}{H_{\mathrm{ctn}}}
\renewcommand{\hom}{\mathrm{Hom}}
\newcommand{\Ialpha}{I_{\alpha}}
\newcommand{\ie}{{\it i.e.,\/}\ }

\newcommand{\limdir}
     {{\displaystyle{\mathop{\rm
     lim}_{\longrightarrow}}}\,}

\newcommand{\Malpha}{M^{\alpha}}
\newcommand{\modA}{\mbox{{\bf mod}-$A$}}
\newcommand{\Nalpha}{N^{\alpha}}
\newcommand{\otimeshat}{\widehat{\otimes}}
\newcommand{\prof}{profinite \relax}
\newcommand{\prf}{\noindent {\bf Proof.\ }}
\newcommand{\qed}{\hfill \ensuremath{\blacksquare} \medskip}

\newtheorem{thm}[subsection]{Theorem}
\newtheorem{prop}[subsection]{Proposition}
\newtheorem{lem}[subsection]{Lemma}
\newtheorem{cor}[subsection]{Corollary}

{\theorembodyfont{\rmfamily}
   \newtheorem{defn}[subsection]{Definition}
   \newtheorem{rem}[subsection]{Remark}
   \newtheorem{example}[subsection]{Example}}

\begin{document}

\begin{center}
   {\Large \bf Cotensor products of modules}\\
   \medskip
   \medskip
   {\large L. Abrams and 
           C. Weibel\footnote{Partially supported by NSF grants.}}\\
   \medskip
   \medskip
   {\small Department of Mathematics, \\
           Rutgers University \\ 
           New Brunswick, NJ 08903}
\end{center}

\medskip

\begin{abstract}
Let $C$ be a coalgebra over a field $k$ and $A$ its dual algebra. The
category of $C$-comodules is equivalent to a category of
$A$-modules. We use this to interpret the cotensor product $M \square
N$ of two comodules in terms of the appropriate Hochschild cohomology
of the $A$-bimodule $M \otimes N$, when $A$ is finite-dimensional,
profinite, graded or differential-graded. The main applications are to
Galois cohomology, comodules over the Steenrod algebra, and the
homology of induced fibrations.
\end{abstract}

\begin{section}{Introduction}

Let $C$ be a coalgebra over a field $k$ and $A$ its dual algebra. It
is classical \cite{Car} that one can identify comodules over $C$ with
a subcategory of modules over $A$. Thus comodule constructions may be
translated into module constructions. In this paper, we interpret the
cotensor product $M \square N$ of two comodules as the vector space
$\hom_{\Ae}(A, M \otimes N)$ of bimodule maps from $A$ to $M \otimes
N$. We also interpret the right derived functors $\cotor_C^*(M,N)$ in
terms of the Hochschild cohomology of the $A$-bimodule $M \otimes
N$. Our results tie together several apparently disparate areas of
mathematics.

When $C$ and $A$ are finite-dimensional, we do indeed have
$\cotor_C^n(M,N) = H^n(A,M \otimes N)$. This straightforward result is
presented in Theorem \ref{th:cotorhoch}.

When $C$ is infinite-dimensional, it is important to think of $A$ as a
topological algebra, as first suggested by Radford \cite{Rad}. In this
case, we prove in Theorem \ref{th:contcotorhoch} that
$\cotor_C^n(M,N)$ is the {\it continuous} Hochschild cohomology
$\Hctn^n(A,M \otimes N)$.

This result has an interesting application to the Galois cohomology of
a profinite group $G$. The profinite group ring $k[G]$ is dual to the
coalgebra $C$ of locally constant functions from $G$ to $k$, and
the Galois cohomology $\Hctn^*(G,M)$ of $G$ with coefficients in
a discrete $G$-module $M$ is the same as $\cotor_C^*(M,k)$;
see Section \ref{sec:profinite} for details.

When $C$ is graded and each $C_n$ is finite-dimensional, we prove that
$\cotor_C^*(M,N)$ is the {\it graded} Hochschild cohomology
$H_{gr}^*(M,N)$, at least when the graded modules $M$ and $N$ are
bounded above (Theorem \ref{th:grcotorhoch}).  Here there is extra
bookkeeping for indices, but no topology is needed.

For example, consider the Steenrod algebra $A^*$ over
$\mathbb{F}_p$, and let $X,Y$ denote topological spaces. Their
homology are comodules over the dual $A_*$, and we have 
$\cotor_{A_*}^n\left(H_*(X),H_*(Y)\right) 
      = H^n\left(A^*,H_*(X \times Y)\right)$; 
see Section \ref{sec:graded}. 

In Section \ref{sec:DG} we extend our results to the differential
graded case, and tie them in with the original formulation and
application of Eilenberg and Moore \cite{EM}.

The special case when $A$ is a Frobenius algebra (and therefore
finite-dimensional) was previously considered by the first author in
\cite{Abr}. The construction there makes use of the special duality
characterizing Frobenius algebras to produce a coalgebra
structure on the dual of $A$ and to translate the comodule
constructions to module constructions. Proposition \ref{pr:injotimes}
here fills in a gap in the proof of \cite[4.6.1]{Abr}.

\smallskip
\noindent
{\bf Notation}

The symbol $\otimes$ denotes tensor product over $k$, and the term
``algebra'' always refers to a $k$-algebra. For any vector space $V$,
write $DV$ for the dual space $\hom_k(V,k)$.  When $V$ is
finite-dimensional there are natural isomorphisms $D(V_1
\otimes V_2) \cong (DV_1) \otimes (DV_2)$ and $V \cong D(DV)$.

\end{section}


\begin{section}{Modules and Comodules}

The results in this section are well known; they appear in \cite{Car}
and later on in \cite{Sw}.  Let $k$ be a field and $A$ a
finite-dimensional algebra.  It is elementary (see \cite[1.1.2]{Sw})
that the dual of the structure maps $k \rightarrow A$ and $A \otimes A
\rightarrow A$ make the dual $DA = \hom_k(A,k)$ into a co-associative
coalgebra.

If $M$ is a left $A$-module, the formula $\Delta_M(m)(a) = am$ defines
a map
\[
\Delta_M \colon M \rightarrow \hom_k(A,M) \cong  M \otimes DA \, .
\]
By \cite[4-03]{Car} or \cite[2.1.2]{Sw}, 
this makes $M$ into a right $DA$-comodule. 
Conversely, if $M$ is a right $DA$-comodule the map 
\[ 
A \otimes M \stackrel{1 \otimes \Delta}{\longrightarrow}
A \otimes (M \otimes DA) \cong (A \otimes DA) \otimes M 
\stackrel{ev \otimes 1}{\longrightarrow}M
\]
makes $M$ into a left $A$-module by \cite[2.1.1]{Sw}.
As observed in {\it loc.~cit.}, these processes are inverse to each
other, and module maps correspond to comodule maps by
\cite{Car} or \cite[2.1.3(e)]{Sw}. We summarize this as follows.

\begin{thm}
   \label{th:modcomod}
Let $A$ be a finite-dimensional algebra. 
Then there is an equivalence
between the category of left $A$-modules and the category of right
$DA$-comodules.
\[
\Amod \cong \comodDA
\]
\end{thm}

Note in particular that the comodule structure map 
$$ \Delta_A \colon A \rightarrow \hom_k(A,A) \cong \enmorph_k(A)$$
is the right regular representation: $\Delta_A(b)(a)=ab$.  
In particular, the composition with $k \rightarrow A$ yields the usual
trace map $\tau_A \colon k \rightarrow A \otimes DA \cong
\enmorph_k(A)$.

When $M$ is finite-dimensional, we can use duality to give another
interpretation of this equivalence.  
If $M$ is any left $A$-module, $\hom_k(M,k)$ is naturally a right
$A$-module by the rule $(fa)(m) = f(am)$.  
If $M$ is finite-dimensional, the dual of the module structure map
$\rho \colon DM \otimes A \rightarrow DM$ is a map $D\rho \colon M
\rightarrow M \otimes DA$; it is straightforward that $D\rho$ makes
$M$ into a right $DA$-comodule.  
Elementary considerations show that this duality $M
\mapsto DM$ gives an equivalence between the category of
finite-dimensional left $A$-modules and the category of
finite-dimensional right $DA$-comodules.

\begin{lem}
When $M$ is finite-dimensional, the structure map $\Delta_M$ agrees
with the duality map $D\rho \colon M\rightarrow M \otimes DA$.
\end{lem}
\prf
Let $\{e_i\}$ be a basis for $A$ and $\{e^i\}$ the dual basis of $DA$.
By inspection, when $M=A$ we have
\[
(D\rho)(1)(e^i \otimes e_j) = \rho(e^i \otimes e_j)(1) = \delta_{ij},
\]
so $(D\rho)(1) = \sum e_i \otimes e^i = \Delta_A(1)$. 
Now given $m \in M$, let $f_m \colon A \rightarrow M$ be the module
homomorphism with $f_m(1)=m$.
By naturality of $\rho$ and $\Delta$,

\medskip\medskip
\begin{minipage}[t]{.9\linewidth}
\begin{center}
$
\Delta_M(m) = ( f_m \otimes DA)\Delta_A(1) = (f_m \otimes
DA)(D\rho)(1) = D\rho(m).
$
\end{center}
\end{minipage}~\qed
\medskip
\medskip

Note that even when $M$ is not finite-dimensional, and hence $D\rho$
might not be a comodule structure map, the map $\Delta_M$ can still be
expressed by $\Delta_M(m) = (f_m \otimes DA)(D\rho)(1)$, where
$f_m\colon A\to M$ is as in the proof of the lemma.

\bigskip
Any construction on $A$-modules translates into a corresponding
construction on $DA$-comodules. The contragradient representation
is one example.

\begin{defn} \label{contragr}
Let $C$ be a finite-dimensional coalgebra. If $M$ is a right
$C$-comodule, $DM=\hom_k(M,k)$ is the left $C$-comodule with 
structure map
\[
\Delta=\Delta_{DM}\colon DM\to DM\otimes C \cong \hom(M,C)
\]
defined by the formula $(\Delta f)(m)=(f\otimes1)(\Delta m)$.
\end{defn}

We leave it to the reader to check that this is the comodule structure
obtained by regarding $M$ as a left module over $A=DC$, and
translating the right $A$-module structure on $DM$ into a left
$DA$-comodule structure. 

\begin{rem} 
   \label{rem:dualmap}
As an exercise, the reader might enjoy verifying that $\Delta_{DM}$ is
the dual of the map $\mu_M\colon A\otimes M\to M$.
\end{rem}

\end{section}

\begin{section}{Cotensor product and Cotor}
   \label{sec:cotor}

In this section, $A$ will be a finite-dimensional $k$-algebra.
If $M$ is a left $A$-module, and $N$ a right $A$-module, then the
results of the last section allow us to view $M$ as a right
$DA$-comodule, and $N$ as a left $DA$-comodule.
As such, we can form the cotensor product $M \square N = M
\square_{DA} N$.
Recall from \cite{EM} that $M \square N$ is defined to be the kernel
of the map 
\[
\phi = \phi_{M,N}\colon M \otimes N \rightarrow M \otimes DA \otimes
N, \ \phi(m \otimes n) = \Delta_1(m) \otimes n - m \otimes
\Delta_2(n).
\]
Here $\Delta_1 \colon M \rightarrow M \otimes DA$ and $\Delta_2 \colon
N \rightarrow DA \otimes N$ are the comodule structure maps.

It is useful to generalize the $\square$ construction to
$A$-bimodules, noting that $M \otimes N$ is an $A$-bimodule. 
If $B$ is an $A$-bimodule, the left and right $A$-module structures
yield maps $\Delta_1 \colon B \rightarrow B \otimes DA$ and $\Delta_2
\colon B \rightarrow DA \otimes B$. Writing $T$ for the twisting
operator, we set 
\[
\phi_B = \Delta_1 - T\Delta_2 \colon B \rightarrow B \otimes DA,
\]
and define $\square(B) := \mathrm{ker}(\phi_B).$ If $B = M \otimes N$,
where $M$ and $N$ are left and right $A$-modules respectively, then
this $\phi_B$ agrees with the previous $\phi_{M,N}$ modulo the natural
twist isomorphism $M \otimes DA \otimes N \cong M \otimes N \otimes
DA$. By definition, $\square(M \otimes N) = M \square N$. 

\begin{prop}   \label{pr:cotensorhom}
If $M$ and $N$ are left and right $A$-modules, respectively, then 
there is a natural isomorphism
\[
M \square_{DA} N  \, \cong \, \hom_{\Ae}(A, M \otimes N)
\]
\end{prop}
\prf
We view $B = M \otimes N$ as a bimodule. Choose a $k$-basis $\{e_i\}$
for $A$ and dual basis $\{e^i\}$ for $DA$. In the particular case $B =
A \otimes A$, the structure maps $\Delta_1, \Delta_2$ are essentially
the right and left regular representations, respectively.  In terms of
the bases we have, for any $b \in B = A \otimes A$,
\[
   \phi_B(b) \ = \ 
   \sum_{i,j} \left( e_ib \otimes e^i \, - \, be_j \otimes e^j \right)
   \ = \
   \sum_i (e_ib - be_i) \otimes e^i
\]
By naturality of $\phi$ we can use a bimodule map $A \otimes A
\rightarrow B$ to show that this formula holds for any bimodule $B$.

If $b \in \square(B)$ then for each $a \in A$ we have $ab=ba$. Hence
there is a bimodule map $f \colon A \rightarrow B$ with $f(a) =
ab$. Since $f(1)=b$, this gives an injection of $\square(B)$ into
$\hom_{\Ae}(A,B)$. 

Conversely, every bimodule map $f \colon A \rightarrow B$ is
determined by $b = f(1)$, and clearly $ab = ba$ for all $a \in A$. For
such $b$, we have $\phi_B(b) = \sum(e_ib -be_i) \otimes e^i = 0$, so
$b \in \square(B)$. Hence every element of $\hom_{\Ae}(A,B) =
\hom_{\Ae}(A, M \otimes N)$ arises in this way.
\qed

To describe the right derived functors of $\square$, we require a
preliminary result on injective modules over finite-dimensional
algebras. 

\begin{prop}   \label{pr:injotimes}
If $E_1$ is an injective left $A$-module, and $E_2$ is an injective
right $A$-module, then $E_1 \otimes E_2$ is an injective left
$\Ae$-module. 
\end{prop}

\prf
Because $A$ is noetherian, every injective module is a direct sum of
indecomposable modules \cite[25.6]{AF}.  Since direct sums of
injective $\Ae$-modules are injective, we may assume that $E_1$ and
$E_2$ are indecomposable injectives, \ie $E_1 = E(A/I_1)$ and $E_2 =
E(A/I_2)$ for some one-sided ideals $I_1$ and $I_2$ (see
\cite[18.12.3]{AF}).  But there are surjections $A^n \rightarrow
D(A/I_i)$, yielding injections $A/I_i \rightarrow (DA)^n$. Hence $E_1$
and $E_2$ are direct summands of some $(DA)^n$. A fortiori, $E_1
\otimes E_2$ is a direct summand of $(DA)^n \otimes (DA)^n$, which is
a sum of copies of the injective $\Ae$-module $DA \otimes DA$.
\qed

\medskip
A hint to the following result appears in \cite[3.1]{Doi}, but the
discussion there remains in the category of comodules.

\begin{thm}   \label{th:cotorhoch}
Let $A$ be a finite-dimensional algebra.
If $M$ is a left $A$-module, and $N$ a right $A$-module, there is a
natural isomorphism
\[
\cotor^*_{DA}(M,N) \cong H^*(A,M \otimes N).
\]
\end{thm}

\prf
This is a formal consequence of \ref{pr:cotensorhom} and
\ref{pr:injotimes}. 
Choosing an injective resolution $M \rightarrow E_1^*$ in $\Amod$, and
$N \rightarrow E_2^*$ in $\modA$, the Cotor groups are the cohomology
of the chain complex
\[
E_1^* \square E_2^* \cong \hom_{\Ae}(A, E_1^* \otimes_k E_2^*).
\]
On the other hand, $M \otimes N \rightarrow E_1^* \otimes E_2^*$ is an
injective resolution in \Aemod, so by \ref{pr:injotimes} the same
complex computes the Hochschild cohomology groups $H^*(A,M \otimes N)
=
\mathrm{Ext}_{\Ae}^*(A,M \otimes N)$.
\qed

\end{section}

\begin{section}{Profinite algebras}
   \label{sec:profinite}

We now generalize the results of \S \ref{sec:cotor} to comodules over
an arbitrary coalgebra $C$. Given any right $C$-comodule $M$,
we can make $M$ into a left module over $A=DC$, via the composite
\cite{Car}\cite[2.1.1]{Sw}
\[
A \otimes M   \stackrel{1 \otimes \Delta}{\longrightarrow}
A \otimes (M \otimes C) \cong
M \otimes (A \otimes C) \stackrel{1 \otimes ev}{\longrightarrow}M.
\]

Similarly, we can make any left $C$-comodule $N$ into a right
$A$-module. As before, we can consider $M \otimes N$ as an
$A$-bimodule. The proof of Proposition \ref{pr:cotensorhom} readily
extends to the current context.

\begin{prop} \label{prop:key}
Let $M$ and $N$ be right and left $C$-comodules, respectively. There
is a natural isomorphism
\[
M \square_C N \cong \hom_{\Ae}(A, M \otimes N)
\]
\end{prop}

The analogue of Theorem \ref{th:cotorhoch} for the Cotor groups
requires a more robust module context, which begins with
the correspondence between arbitrary coalgebras and
profinite algebras.

Any coalgebra $C$ is the union of its finite-dimensional subcoalgebras
$\Calpha$, so its dual $A = \hom_k(C,k)$ is the inverse limit of
the finite-dimensional algebras $\Aalpha = D(\Calpha)$. We view $A =
\{ \Aalpha \}$ as a pro-object in the category of finite-dimensional
algebras, a structure we will refer to as a {\it \prof algebra.} Each
$\Aalpha$ is isomorphic to $A/\Ialpha$, where $\Ialpha$ is the ideal
of functions vanishing on $\Calpha$; these ideals define a topology on $A$.

Conversely, suppose that $A = \{ \Aalpha \}$ is a \prof algebra.
Then the union $C$ of the filtered system of duals $\Calpha=D(\Aalpha)$ 
is a coalgebra, and clearly $A$ is the \prof algebra associated to $C$. 
This establishes an equivalence between the category of coalgebras and
the opposite category of \prof algebras \cite{Br}\cite{Wit}.

Note that a \prof algebra carries more structure,
(\eg topology) than its inverse limit. This is clearly seen in the
papers by Radford and Witkowski on reflexive coalgebras
\cite{Rad}\cite{Wit}.

\begin{defn}
Let $A = \{ \Aalpha \}$ be a \prof algebra.  We say that a (left)
$A$-module $M$ is {\it rational} if each element of $M$ generates a
finite-dimensional $A$-submodule isomorphic to a quotient of some
$\Aalpha$. 
\end{defn}

The category $\Amodrat$ of rational modules and
$A$-module maps is an abelian category, and the inclusion $\Amodrat
\subset \Amod$ is exact. The argument of \cite[6.11.10]{Wei} shows
that $\Amodrat$ has enough injectives.

Radford \cite[2.2]{Rad} and Witkowski \cite{Wit} have shown that our
notion of rational module agrees with Sweedler's notion \cite[p. 37]{Sw}. 
That is, an $A$-module $M$ is rational in our sense if and only if the
map $M\stackrel{\rho}{\to} \hom(A,M)$ defined by $\rho(m)(a)=a\cdot m$
lands in the subspace $M\otimes C$ of $\hom(A,M)$. 

If $M$ is a rational module, we can regard it as a $C$-comodule as
follows. Since $M=\bigcup Am$, it suffices to consider $Am$. Choose
$\alpha$ so that $Am$ is an $\Aalpha$-module. 
By Theorem \ref{th:modcomod}, $Am$ is a comodule over $C^\alpha$ and
hence over $C$. This proves Sweedler's theorem \cite[2.1.3]{Sw}, which
we record here.

\begin{thm}   \label{th:modcomodprofinite} 
If $C$ is a coalgebra and $A$ the dual
\prof algebra, then there is an equivalence between the category of
right $C$-comodules and the category of left rational $A$-modules.
\[
\Amodrat \cong \mbox{{\bf comod}-}C
\]
\end{thm}

\begin{example}
Let $G$ be a profinite group, with finite quotients $\Galpha$.  The
group rings $k[\Galpha]$ are dual to the Hopf algebras $k^{\Galpha}$
of functions $\Galpha \rightarrow k$.  Hence the \prof
group algebra $\{k[\Galpha]\}$ corresponds to the coalgebra $C =
\bigcup k^{\Galpha}$ of locally constant functions $G \rightarrow k$.
Regarding each $\Galpha$ as the algebraic group
$\mathrm{Spec}(k^{\Galpha})$, we get a pro-algebraic group 
$\{\Galpha\}$.

A rational $k[G]$-module is the same thing as a 
{\it discrete $G$-module} in the sense of Galois cohomology \cite{Wei}. 
By Sweedler's theorem, it is also just a comodule for $C$.
Now a {\it rational representation} of the pro-algebraic group
$\{\Galpha\}$ is just a union of rational representations of the
$\Galpha$, each of which is just a $k^{\Galpha}$-comodule 
(see \cite[2.23]{Fog}). 
Thus a discrete $G$-module may be thought of
as a rational representation of $\{\Galpha\}$.
\end{example}

\begin{example}\label{ex:locfin} (Taft \cite{T}) An $A$-module $M$ is
called {\it locally finite} if each element of $M$ generates a
finite-dimensional $A$-submodule. This notion does not involve
the topology on $A$, and is weaker than the notion of rational
module. 

For example, let $V$ be an infinite-dimensional $k$-vector space and
form the coalgebra $C := k \oplus V$, where the elements of $V$ are
primitive. The dual algebra is $A=k \oplus V^*$. The $A$-module
$M:=V^{**}$ is locally finite (for each $m \in M$, the submodule $Am$
is two-dimensional) but not rational ($Am$ 
need not be a quotient of any $\Aalpha$).
\end{example}


We now turn to $A$-bimodules. In order to define continuous Hochschild
cohomology, we again need to consider the topology of $A$.

\begin{defn}
We say that an $A$-bimodule $B$ is {\it rational} if each
element of $B$ generates a finite-dimensional sub-bimodule isomorphic
to a quotient of some $\Aalpha$. That is,
$B$ is a rational module over the \prof algebra
$A \otimeshat A^{\mathrm{op}} = 
\{ \Aalpha \otimes \Aalpha^{\mathrm{op}} \}$.

We define the {\it continuous Hochschild cohomology} of $B$ to be the
right-derived functors $\Hctn^*(A,B)$ of
\[
B \mapsto \hom_{\Ae}(A,B) = 
 \left\{b \in B : ab = ba \ \mbox{for all}\ a \in A \right\}.
\]
\end{defn}

Recall \cite[6.5.1]{Wei} that if $B$ is an $A$-bimodule, its
Hochschild cohomology $H^*(A,B)$ is defined to be the homology of the
cochain complex $C^*(A,B)$, where $C^n(A,B)$ denotes the $k$-module of
$(n+1)$-fold multilinear maps from $A$ to $B$, \ie maps $f \colon
A^{\otimes n} \rightarrow B$. When $n=0$, $C^0(A,B) = B$.

When $A = \{ \Aalpha \}$ is a \prof algebra and $B$ is a rational
bimodule, we define $\Cctn^n(A,B)$ to be the subspace of $C^n(A,B)$
consisting of {\it continuous cochains}, \ie maps which factor through
some quotient $\Aalpha^{\otimes n}$ of $A^{\otimes n}$.

For each finite quotient $\Aalpha = A/\Ialpha$ of $A$, set $\Balpha :=
\{ b \in B : \Ialpha b = b \Ialpha = 0 \}$. Then $\Balpha$ is an
$\Aalpha$-bimodule, and $B = \bigcup \Balpha$ if $B$ is rational. 

\begin{prop}
   \label{pr:respectlim}
We have $\Cctn^*(A,B) \ = \ \limdir C^*(\Aalpha,\Balpha)$.
\end{prop}
\prf
The composite of $A^{\otimes n} \rightarrow \Aalpha^{\otimes n}$ with
a map $\Aalpha^{\otimes n} \rightarrow \Balpha$ is a continuous
cochain, so $C^*(\Aalpha,\Balpha) \rightarrow C^*(A,B)$ is an
injection. Conversely, given a continuous cochain $f \colon A^{\otimes
n} \rightarrow B$, there is an $\alpha$ so that $f$ factors through
$\Aalpha^{\otimes n}$. Because $\Aalpha^{\otimes n}$ is
finite-dimensional and $B$ is rational, the image $f(A^{\otimes n})$
lies in some $B^{\beta}$. Choosing $\gamma$ so that $A_{\gamma}$ maps
to $\Aalpha$ and $A_{\beta}$, we see that $f \in
C^n(A_{\gamma},B^{\gamma})$.
\qed

\begin{rem}
Our terminology comes from the fact that the \prof algebra $A
\otimeshat \cdots \otimeshat A = \{ \Aalpha^{\otimes n} \}$ may be
regarded as a topological algebra. If $B$ has the discrete topology, a
continuous map $A \otimeshat \cdots \otimeshat A \rightarrow B$ must
factor through some $\Aalpha^{\otimes n}$ and hence be a continuous
cochain.
\end{rem}

\begin{thm}
   \label{th:ctnlimit} 
If $B$ is a rational bimodule over a \prof algebra $A$, then
$\Hctn^*(A,B)$ is the cohomology of the complex $\Cctn^*(A,B)$. In
particular,
\[
\Hctn^*(A,B) \cong \limdir H^*(\Aalpha,\Balpha).
\]
\end{thm}
\prf
For simplicity, set 
\[
T^n(B) = H^nC^*(A,B) = \limdir H^nC^*(\Aalpha,\Balpha).
\]
Note that $C^0(A,B) \cong B$, and thus $T^0(B) = \Hctn^0(A,B)$. We now
argue as in the proof of \cite[6.11.13]{Wei}. The set $\{T^n\}$ forms
a $\delta$-functor because $H^*C^*(\Aalpha,\Balpha)$ is the Hochschild
cohomology $H^*(\Aalpha,\Balpha)$ for each $\alpha$. To see that this
$\delta$-functor is universal, note that if $J$ is an injective object
in $A \otimeshat \Amodrat$ then each $J^{\alpha}$ is an injective
$\Aalpha$-bimodule (because $B \mapsto B^{\alpha}$ is right adjoint to
the forgetful functor). Hence if $n \neq 0$ then
\[
\displaystyle{T^n(J) = \limdir H^n(\Aalpha,J^{\alpha}) = 0}.
\hspace{3em} \qed
\]

\begin{example}
Let $G = \{ G_{\alpha} \}$ be a profinite group and $M$ a discrete
left  $G$-module. A continuous cochain $f \colon G^n \rightarrow M$
in the sense of Galois cohomology is one which factors through some
quotient $G^n_{\alpha}$ of $G^n$. Its linear extension $kG^{\otimes n}
\rightarrow M_{\epsilon}$ is a continuous cochain, where
$M_{\epsilon}$ is the bimodule with trivial right $G$-action. Thus
$\Cctn^*(kG,M_{\epsilon})$ is the chain complex used to compute the
Galois cohomology of $G$. It follows that 
$\Hctn^*(k[G],M) \cong \Hctn^*(G,M)$.
\end{example}

We can now prove the main theorem.

\begin{thm}
   \label{th:contcotorhoch}
If $C$ is a coalgebra, $A$ its \prof dual, and $M,N$
are $C$-comodules, then
\[
\cotor_C^*(M,N) = \Hctn^*(A, M \otimes N) .
\]
\end{thm}
\prf
Let $\{ \Calpha \}$ and $\{ \Aalpha \}$ be as in the beginning of this
section. Define $\Malpha := \{ x \in M \mid \Ialpha x = 0 \} = \{x \in
M \mid \Delta_Mx \in M \otimes \Calpha\}$ and define $\Nalpha$
similarly. Note that $\Malpha \otimes \Nalpha = (M \otimes
N)^{\alpha}$. By Theorem \ref{th:modcomodprofinite} $M, N, \Malpha,
\Nalpha$ may all be viewed as rational $A$-modules. 

Note that $\cotor_{\Calpha}(\Malpha, \Nalpha)$ can be calculated with
the complex
\[
\Malpha \otimes \Nalpha 
\rightarrow \Malpha \otimes \Calpha \otimes \Nalpha 
\rightarrow \Malpha \otimes \Calpha \otimes \Calpha \otimes \Nalpha 
\rightarrow \cdots
\]
Taking direct limits, we get a complex for $\cotor_C(M,N)$. Since
homology commutes with this particular limit, we have 
\[
\cotor_C^*(M,N) \ = \ \limdir \cotor_{\Calpha}^*(\Malpha,\Nalpha)
\]
Applying Proposition \ref{pr:respectlim} and 
Theorem~\ref{th:cotorhoch} gives
\[
 \limdir \cotor_{\Calpha}^*(\Malpha,\Nalpha)
   \ = \  \limdir H^*(\Aalpha,\Malpha \otimes \Nalpha)
   \ = \  \Hctn^*(A;M \otimes N),
\]
which completes the proof. \hfill \qed

\end{section}

\begin{section}{Graded modules and comodules}
   \label{sec:graded}

Replacing the category of vector spaces with the category of graded
vector spaces does not change things very much, as we now explain.

First, we need to fix our notation, for which we follow \cite{Mac} and
\cite{EM}. 
The category of graded vector spaces has an internal $\hom$,
constructed as follows.  Given two graded vector spaces $V$ and $V'$,
the degree $p$ component $\hom^p_k(V,V')$ of $\hom^*_k(V,V')$ is the
vector space of all maps $f\colon V\to V'$ of degree $p$, i.e.,
satisfying $f(V^i)\subset (V')^{i+p}$. The indexing is set up so that
the evaluation map $\hom_k^*(V,V')\otimes V \to V'$ is homogeneous of
degree zero.  The {\it graded dual} $DV=\hom^*(V,k)$ of $V$ arises as
the special case $V'=k$. Thus $(DV)_n=(DV)^{-n}=\hom_k(V^n,k)$.  If
each $V^n$ is finite-dimensional we say that $V$ is of {\it finite
type}; this is the hypothesis needed to have $D(DV)\cong V$.

By a graded algebra $A$ we will mean a positively graded algebra
$A^0\oplus A^1\oplus\cdots$. 
We will assume throughout this section that $A$ is of finite type.

By a graded coalgebra $C$ we will mean a positively
graded coalgebra $C_0\oplus C_1\oplus\cdots$.
If $A$ is a graded $k$-algebra of finite type, then
$DA$ is a co-associative graded coalgebra by
\cite[3.1(4)]{MM} or \cite[6.0.2]{Sw}.
This is because the tensor product $A\otimes A$ is of finite type.

Note that a graded module $M$ is rational if and only if it is
locally finite (example~\ref{ex:locfin}); 
the topology on $A$ plays no role when $A$ is of finite type. 
As a typical example, if $M$ is a
{\it bounded-above} module (\ie $M^p = 0 \mbox{\ for\ } p \gg 0$) then
clearly $M$ is locally finite, because $A$ has finite type.

\begin{prop}   \label{pr:grmodcomod} 
If $A$ is of finite type, then the category of graded right
$DA$-comodules is equivalent to the category of locally finite, graded
left $A$-modules. Under this equivalence, the category of
bounded-below comodules $(M_p=0 \mbox{\ for \ } p \ll 0)$ is
equivalent to the category of bounded-above modules.
\end{prop}

\prf
If $M_*$ is a comodule, the classical formula for the $A$-module
structure on $M$ makes it a locally finite graded module
(see \cite{MM}). Conversely, if $M^*$ is a
graded module, the formula $\Delta_M(m)(a) = am$ is homogenous, and
$\Delta_M(m)$ is in $M \otimes DA$ if $M$ is locally finite.
\qed

Suppose now that $M$ is a graded right comodule and $N$ a graded left
comodule. Then we can consider $M$ and $N$ as left and right
modules, respectively, and form the graded bimodule $M\otimes N$. 
Applying the graded $\hom_{\Ae}(A,-)$ (defined as in 
\cite[p.~185]{Mac}) allows us to construct the graded Hochschild
cohomology $H_{gr}^*(A,M\otimes N)$, as in \cite[p.~300]{Mac}.

The proof of Proposition \ref{pr:cotensorhom} goes through to prove
the following analogue. 

\begin{prop} \label{prop:graded} Let $A$ be a graded algebra of finite
type, with dual coalgebra $C$.
If $M$ and $N$ are locally finite left and right graded $A$-modules,
respectively, then there is a natural isomorphism
\[
M \square_{C} N  \, \cong \, \hom_{\Ae}(A, M \otimes N)
\]
\end{prop}

In order to prove the analogue of Theorem \ref{th:cotorhoch}, we
require some preliminary results.

\begin{lem}\label{lem:adjoint}
If $V$ is a bounded above graded vector space, then $DA \otimes V$ is
an injective graded $A$-module (and hence also an injective graded
$DA$-comodule), and $DA \otimes V \otimes DA$ is an
injective graded $A$-bimodule.
\end{lem}
\prf
Because $\hom_k(A,-)$ is right adjoint to the forgetful functor from
graded $\!A$-modules to graded 
vector spaces, it preserves injectives \cite[2.3.10]{Wei}.
 Hence $\hom_k(A,V)$ is injective. In general, 
$\hom^n_k(A,V) = \prod_{i+j=n} DA^i \otimes V^j$ is larger than 
$(DA\otimes V)^n = \bigoplus DA^i \otimes V^j$. 
But if $V$ is bounded above, we have $DA \otimes V = \hom(A,V)$, which
is injective, and similarly
\[
DA \otimes V \otimes DA = D(A \otimes A) \otimes V = \hom_k(\Ae,V)
\]
is injective as a graded $\Ae$-module, \ie a graded bimodule.
\qed

\begin{rem}
If $A$ is a Noetherian algebra, then direct sums of injective modules
are injective. Since $DA$ and $DA \otimes DA$ are injective (compare
with the proof of Proposition \ref{pr:injotimes}), we may
remove the hypothesis that $V$ be bounded above.
\end{rem}

\begin{cor}   \label{cor:res} 
If $M$ and $N$ are bounded-above graded $A$-modules, there are
bounded-above graded vector spaces $V^i$ and $W^i$ and injective
resolutions $M \rightarrow I^*, \, N \rightarrow J^*$ with $I^i \cong
DA \otimes V^i, \, J^i \cong W^j \otimes DA$. Moreover, $M \otimes N
\rightarrow I^* \otimes J^*$ is an injective bimodule resolution.
\end{cor} 

\begin{thm}   \label{th:grcotorhoch} 
Suppose $A$ is a graded $k$-algebra of finite type, with dual
coalgebra $C$, and that $M$, $N$ are bounded-above graded
$A$-modules. Then for each $n$ there is a natural 
isomorphism of graded vector spaces:
\[
\cotor_C^n(M,N) \cong H_{gr}^{n}(A, M \otimes N).
\]
\end{thm}
\prf
Take injective resolutions $M \rightarrow I^*$ and $N \rightarrow J^*$
as in Corollary \ref{cor:res}. These lie in the category of graded
$DA$-comodules. Because $I^p$ and $J^q$ are bounded below, the proof
of Proposition \ref{pr:cotensorhom} shows that $I^p \square J^q =
\hom_{\Ae}(A, I^p \otimes J^q)$. Thus we have
\begin{eqnarray*}
\cotor_C^n(M,N) & = & H^n(I^* \square J^*) \\
                & = & H^n\hom_{\Ae}(A,I^* \otimes J^*) \\
                & = & H_{gr}^n(A, M \otimes N) \ \ \ \  \qed
\end{eqnarray*}

\begin{rem} Eilenberg and Moore actually use homological indexing,
defining $\cotor^C_n(M,N)$ to be the product (over $p+q=n$)
of the $\cotor_{p,q}(M,N)$, which in our notation is the
homological degree~$p$ part of the graded vector space
$\cotor_C^{-q}(M,N)$; see \cite[p.~207]{EM}.
\end{rem}

\begin{example}
Let $A^*$ denote the Steenrod algebra over $\mathbb{F}_p$. If $X$ is a
topological space, its cohomology $H^*(X)$ is a bounded-below module
over $A^*$ and its homology $H_*(X)$ is a graded comodule over the
dual $A_*$. We can view $H_*(X)$ as cohomologically bounded above, and
hence by \ref{pr:grmodcomod} as a bounded-above module over $A^*$,
with $H^*(X) = DH_*(X)$. Because $H_*(X \times Y) \cong H_*(X)
\otimes H_*(Y)$ is a bicomodule over $A^*$, we have
\[
H_*(X) \square_{A_*} H_*(Y) 
    =    \hom_{\Ae}\left(A^*,H_*(X \times Y)\right).
\]
and
\[
\cotor_{A_*}^n\left(H_*(X),H_*(Y)\right) 
      = H^n\left(A^*,H_*(X \times Y)\right).
\]
\end{example}

\end{section}


\begin{section}{The Differential Graded Case} 
   \label{sec:DG}

We now pass to the differential graded setting, using essentially the
same notation as in the graded case.
By a differential graded algebra (or DG-algebra) we will mean a graded
algebra $A$ with a differential of degree one, satisfying the usual
Leibniz relation for $d(aa')$.  
We will regard $A$ as a cochain complex:
\[
0 \rightarrow A^0 \rightarrow A^1 \rightarrow \cdots.
\]
We use cochain notation $M^*$ for DG $A$-modules, so their
differential has degree one.

By a differential graded coalgebra (or DG-coalgebra) $C$ we will mean
a positively graded coalgebra with a differential of degree $-1$,
satisfying the co-Leibniz relation.  We will regard $C$ as a chain
complex:
\[
 \cdots \rightarrow C_1 \rightarrow C_0 \rightarrow 0.
\]
For example, $C=DA$ is a DG-coalgebra with $C_n=\hom(A^n,k)$. We use
chain complex notation for DG $C$-comodules, so their differential has
degree $-1$.

\begin{prop}
If $A$ is of finite type, then the category of right DG $DA$-comodule
is equivalent to the category of locally-finite left DG $A$-modules.
\end{prop}
\prf
The only point which must be added to the proof for the graded case
(Proposition \ref{pr:grmodcomod}) is that the same differential is
used for both $A$-modules and $DA$-comodules.
\qed

If $M$ and $N$ are right and left DG-comodules, then their cotensor
product $M\square_CN$ is naturally a DG-vector subspace of $M\otimes
N$.  Yet again, the proof of Proposition \ref{pr:cotensorhom} goes
through to prove the following DG analogue.

\begin{prop} \label{prop:DGkey} Let $A$ be a DG algebra of finite
type, with dual coalgebra $C$.
If $M$ and $N$ are locally finite left and right DG $A$-modules,
respectively, then there is a natural isomorphism of DG vector spaces:
\[
M \square_{C} N  \, \cong \, \hom_{\Ae}(A, M \otimes N).
\]
\end{prop}

The cotorsion groups $\cotor^C_n(M,N)$ were defined by Eilenberg and
Moore in \cite[p.~206]{EM} as follows.  If $M\to I^*$ and $N\to J^*$
are injective DG comodule resolutions, then $I^*\square_C J^*$ is a
triple chain complex whose $(i,j,k)$ entry is $(I^{-i}\square_C
J^{-j})_{k}$. Here, the index $k$ refers to the total internal grading
of the DG comodules $I^{-i},J^{-j}$. The group $\cotor^C_n(M,N)$ is
defined to be the $n$th homology of the (product) total complex of
$(I^{-i}\square_C J^{-j})_{k}$. The homological indexing reflects not
only the applications to homology of spaces (see
\ref{ex:EMresult} below), but the fact that we want
$\cotor^C_n(C,C)\cong H_n(C)$ to be positively indexed.

If $M$ and $N$ are left and right DG $A$-modules,
their tensor product $M\otimes N$ may be viewed as a DG bimodule.
As such, we may form the differential graded Hochschild cohomology
$H_{DG}^{p}(A,M\otimes N)$.

\begin{thm}   \label{thm:dgcotor}
Suppose $A$ is a DG algebra of finite type with dual coalgebra $C$
and $M$, $N$ are bounded-above DG $A$-modules. Then for each $n$ there
is a natural isomorphism:
\[
\cotor^C_n(M,N) \cong H_{DG}^{-n}(A,M\otimes N).
\]
\end{thm}

\prf
Both Lemma \ref{lem:adjoint} and \ref{cor:res}
carry over to the DG setting with almost no change. Thus we can choose
injective objects $I^*$ and $J^*$ in the category of DG-modules
over $A=DC$ such that $M\otimes N\to I^*\otimes J^*$ is an
injective  resolution in the category of DG-bimodules.

By definition, $H_{DG}^{n}(A,M\otimes N)$ is the $n$th cohomology of
the (product) total complex of the triple cochain complex
$\hom_{A^e}(A,I^*\otimes J^*)$. The $(i,j,k)$ entry of this complex
is
\[
\hom_{A^e}(A,I^i\otimes J^j)^k = (I^i\square_C J^j)_{-k}.
\]
The standard conversion between chain complexes and cochain complexes
($V^n=V_{-n}$) identifies this triple complex with the triple chain
complex used above to define $\cotor^C(M,N)$.
\qed

\begin{example} \label{ex:EMresult} (Eilenberg-Moore)
Suppose given a cartesian square of CW complexes
\begin{eqnarray*}	E' &\to& E \\ 
\downarrow && \downarrow\pi \\
	B' &\to& B
\end{eqnarray*}
in which $B$ is simply connected, with finitely many cells in each
dimension, and $\pi$ is a fibration.  Let $C_*B$ be the cellular
complex of $B$ with coefficients in a field $k$, and
$C^*B=\hom(C_*B,k)$ its dual algebra.  This is a DG algebra of finite
type.  The main result (theorem~12.1) of \cite{EM}, combined with our
Theorem \ref{thm:dgcotor} (and the Eilenberg-Zilber theorem) states
that the homology of $E'$ is isomorphic to
\[
H_*(E';k) \cong \cotor^{C_*B}(C_*E,C_*B') \cong
H^*_{DG}(C^*B, C_*(E\times B')).
\]

In particular, if $F$ is the fiber of the fibration $E\to B$, then
(taking $B'$ to be a point) we have
\[
H_*(F;k) \cong \cotor^{C_*B}(C_*E,k) \cong
H^*_{DG}(C^*B, C_*E\otimes k).
\]
Finally, if $E$ is contractible then $C_*E \cong k$ and 
\[
H_*(\Omega B;k) \cong \cotor^{C_*B}(k,k) \cong H^*_{DG}(C^*B,k) .
\]
\end{example}

\end{section}


\end{document}